\documentclass[11pt]{article}
\usepackage[T1]{fontenc}
\usepackage{latexsym,amssymb,amsmath,amsfonts,amsthm}
\newcommand{\R}{{\mathbb R}}
\newcommand{\Z}{{\mathbb Z}}

\newcommand{\C}{{\mathbb C}}

\newcommand{\ds}{\displaystyle}
\newcommand{\no}{\nonumber}
\newcommand{\be}{\begin{eqnarray}}
\newcommand{\ben}{\begin{eqnarray*}}
\newcommand{\en}{\end{eqnarray}}
\newcommand{\enn}{\end{eqnarray*}}

\newcommand{\bt}{\beta}

\newcommand{\curl}{{\rm curl\,}}
\newcommand{\dive}{{\rm div\,}}
\newcommand{\I}{{\rm Im\,}}

\newcommand{\f}{\frac}
\newcommand{\G}{\Gamma}

\newcommand{\Om}{\Omega}
\newcommand{\om}{\omega}

\newcommand{\al}{\alpha}

\newtheorem{theorem}{Theorem}[section]
\newtheorem{lemma}[theorem]{Lemma}

\newtheorem{remark}[theorem]{Remark}

\begin{document}
\renewcommand{\theequation}{\arabic{section}.\arabic{equation}}
%\begin{titlepage}
\title{\bf The linear sampling method for the inverse electromagnetic scattering
by a partially coated bi-periodic structure}
\author{Guanghui Hu and Bo Zhang\\
LSEC and Institute of Applied Mathematics\\
Academy of Mathematics and Systems Science\\
Chinese Academy of Sciences\\
Beijing 100190, China\\
({\sf huguanghui@amss.ac.cn} (GH),\ {\sf b.zhang@amt.ac.cn} (BZ))}
\date{}
%\end{titlepage}
\maketitle

%\vspace{.2in}
%\noindent
%Proposed running head: {\bf Linear sampling method for inverse scattering by bi-periodic structures}
%\vspace{.1in}
%
%\noindent
%Manuscript Correspondence Address:\\
%Professor Bo Zhang\\
%Institute of Applied Mathematics\\
%Academy of Mathematics and Systems Science\\
%Chinese Academy of Sciences\\
%Beijing 100190, China\\
%Email: b.zhang@amt.ac.cn\\
%Telephone: +86 10 6265 1358\\
%Fax:       +86 10 6254 1689
%
%\newpage
%\begin{center}
%\title{\Large\bf The linear sampling method for the inverse electromagnetic scattering
%by a partially coated bi-periodic structure}
%\end{center}
%
%\vspace{.2in}

\begin{abstract}
In this paper, we consider the inverse problem of recovering a doubly periodic Lipschitz
structure through the measurement of the scattered field above the structure produced
by point sources lying above the structure.
The medium above the structure is assumed to be homogenous and
lossless with a positive dielectric coefficient. Below the structure
is a perfect conductor partially coated with a dielectric.
A periodic version of the linear sampling method is developed to reconstruct the
doubly periodic structure using the near field data.
In this case, the far field equation defined on the unit ball of
$\mathbb{R}^3$ is replaced by the near field equation which is a
linear integral equation of the first kind defined on a plane above the periodic surface.

\vspace{.2in}
\end{abstract}

\section{Introduction}
\setcounter{equation}{0}

Consider the problem of scattering of electromagnetic waves by a doubly periodic structure
of period $\Lambda=(2\pi,2\pi)$ defined by
\ben
\G=\{x_3=f(x_1,x_2)\,|\,f(x_1+2n_1\pi,x_2+2n_2\pi)=f(x_1,x_2)>0,
\quad\forall\;n=(n_1,n_2)\in\Z^2\},
\enn
where the function $f$ is assumed to be Lipschitz continuous so the periodic structure
$\G$ is a Lipschitz surface.
The medium above the structure is assumed to be homogenous with a constant dielectric
coefficient $\epsilon_0>0$, and below the structure is a perfect conductor with a partially
coated dielectric boundary. The magnetic permeability is assumed to be a positive constant
$\mu_0$ throughout $\R^3$. Given the structure and a time-harmonic electromagnetic
wave incident on the structure, the direct scattering problem is to compute
the electric and magnetic distributions away from the structure.
In this paper, we are interested in the inverse problem of reconstructing the shape of
the bi-periodic structure from a knowledge of the incident and scattered fields.
The purpose of this paper is to develop a periodic version of the Linear Sampling Method
for such an inverse problem. We refer to \cite{PR} for historical remarks and details
of the applications of the scattering theory in periodic structures
and \cite{CC06} for a recent overview of the linear sampling method.

Physically, the propagation of time-harmonic electromagnetic waves (with the time variation
of the form $e^{-i\om t},$ $\om>0$) in a homogeneous isotropic medium in $\R^3$ is modeled
by the time-harmonic Maxwell equations:
\be\label{1.1}
 \curl E-ikH=0,\qquad\;\; \curl H+ikE=0.
\en
Here, we assume that the medium above the structure is lossless, that is,
$k$ is a positive wave number given by $k=\sqrt{\epsilon_0\mu_0}\omega$
in terms of the frequency $\omega,$ the electric permittivity $\varepsilon_0$
and the magnetic permeability $\mu_0$.
Consider the time-harmonic plane wave
\ben
E^i=pe^{ikx\cdot d},\qquad  H^i=qe^{ikx\cdot d},
\enn
incident on $\G$ from the top region $\Omega:=\{x\in\R^3\,|\,x_3>f(x_1,x_2)\}$,
where $d=(\al_1,\al_2,-\bt)=(\cos\theta_1\cos\theta_2,\cos\theta_1\sin\theta_2,-\sin\theta_1)$
is the incident direction specified by $\theta_1$ and $\theta_2$ with
$0<\theta_1<\pi,\,0<\theta_2\leq 2\pi$ and the vectors $p$ and $q$ are polarization directions
satisfying that $p=\sqrt{{\mu}/{\varepsilon}}(q\times d)$ and $q\bot d.$

In this paper, we assume that the boundary $\G$ has a Lipschitz dissection
$\G=\G_D\cup\Sigma\cup\G_I,$ where $\G_D$ and $\G_I$ are disjoint,
relatively open subsets of $\G$ having $\Sigma$ as their common boundary.
Suppose below $\G$ is a perfect conductor partially coated by a dielectric on $\G_I$.
The problem of scattering of time-harmonic electromagnetic waves is modeled by
the following exterior mixed boundary value problem:
\be\label{equation1}
  \curl\curl E-k^2E=0  &&\rm{in}\quad\Om,\\ \label{equation2}
  \nu\times E=0         &&\rm{on}\quad\G_D,\\ \label{equation3}
  \nu\times\curl E-i\lambda (\nu\times E)\times\nu=0 &&\rm{on}\quad\G_I,\\ \label{equation4}
   E=E^i+E^s  &&\rm{in}\quad\Om,
\en
where $\nu$ is the unit normal of $\G$ pointing into $\Om$.
We assume throughout this paper that $\lambda$ is a positive constant.

Set $\al=(\al_1,\al_2,0)\in\R^3,\,n=(n_1,n_2)\in\Z^2$.
Motivated by the periodicity of the medium we look for $\al$-quasi-periodic
solutions in the sense that $E(x_1,x_2,x_3)e^{-i\al\cdot x^\prime}$ is
$2\pi$ periodic with respect to $x_1$ and $x_2$, respectively.
Since the domain is unbounded in the $x_3$-direction, a radiation condition
must be imposed. Physically it is required that the scattered fields remain
bounded as $x_3$ tends to $+\infty$, which leads to the so-called outgoing
wave condition of the form:
\be\label{RE}
E^s(x)=\sum_{n\in\Z^2}E_ne^{i(\al_n\cdot x+\bt_nx_3)},\qquad
x_3>\max_{x_1,x_2}f(x_1,x_2),
\en
where $\al_n=(\al_1+n_1,\al_2+n_2,0)\in\R^3$, $E_n=(E_n^{(1)},E_n^{(2)},E_n^{(3)})\in\C^3$
are constant vectors and
$$
\beta_n=\left\{\begin{array}{lll}
             (k^2-|\al_n|^2)^{\frac{1}{2}}\qquad\rm{if}\;|\al_n|<k,\\
              i(|\al_n|^2-k^2)^{\frac{1}{2}}\qquad\rm{if}\;|\al_n|>k,\\
             \end{array}\right.
$$
with $i^2=-1$. Furthermore, we assume that $\beta_n\neq 0$ for all $n\in\Z^2$.
The series expansion in (\ref{RE}) will be considered as the Rayleigh series
of the scattered field, and the condition is called the Rayleigh expansion radiation condition.
The coefficients $E_n$ in (\ref{RE}) are also called the Rayleigh sequence.
From the fact that $\rm{div}\,E^s(x)=0$ in $\Om$ it is clear that
\be\label{div}
   \al_n\cdot E_n+ \beta_nE_n^{(3)}=0.
\en

The inverse problem considered in this paper is concerned with
determining the profile $\G$ and the impedance coefficient $\lambda$
from a knowledge of the incident wave $E^i$ and the tangential component of the total electric
field, $\nu\times E$, on a plane $\G_b=\{x\in\R^3\,|\,x_3=b\}$ above the structure.
The uniqueness of this inverse problem was proved in \cite{HQZ} for the case when the incident
waves are electric dipoles. Precisely, it was shown in \cite{HQZ} that,
if the tangential components on $\G_b$ of two scattered electric fields
are identical for all incident electric dipoles $E^{in}(x;y)=\curl_x\curl_x\{PG(x,y)\}$
with all $y\in\G_b$ and three linear independent vectors $P\in\R^3$,
then their corresponding scattered periodic structures $\G_j(j=1,2)$ and the impedance
coefficients $\lambda_j(j=1,2)$ on $\G_I$ must coincide, where $G(x,y)$ is
the free-space quasi-periodic Green function (see Section \ref{sec2}).
In this paper, we are interested in numerically reconstructing the shape of the periodic
structure $\G$ by using the idea of the linear sampling method.
The linear sampling method was proposed in \cite{CA} for numerically reconstructing
the shape and location of the obstacle in the inverse acoustic obstacle scattering problems.
This method has attracted extensive attention in recent years since it does not need
to know the physical property of the scattering obstacles in advance.
The application of the linear sampling method to the inverse electromagnetic scattering
problems can be found in \cite{CCM,CMH,CHM}.
Recently in \cite{HQZ2}, a periodic version of the linear sampling method was proposed
and implemented for the two-dimensional TE polarization case of the inverse problem
considered in this paper, where the Maxwell equations are replaced by the scalar Helmholtz
equation and the boundary conditions on $\G_D$ and $\G_I$ are replaced with the
Dirichlet and impedance conditions, respectively.
In \cite{Kirsch1998}, Kirsch proposed a mathematically-justified version of the linear
sampling method, the so-called factorization method. However, it is still an open question
to characterize a bounded conducting obstacle for the Maxwell equations by using the
factorization method (see \cite{Kirsch2004}). We refer to \cite{AN,AK,L2010} for the application
of the factorization method to the 2D inverse problems by diffraction gratings with the
Dirichlet, impedance and transmission conditions and to \cite{Tilo99} for a recent convergence
result of the linear sampling method as well as a connection between the linear sampling
and factorization methods.

The inverse scattering problem by a smooth doubly periodic structure has been studied in \cite{AH,BZ}
for the case when $\G_I=\emptyset$. With a lossy medium (i.e., $\I(k)>0$) above the conductor,
Ammari \cite{AH} proved a global uniqueness result for the inverse problem with one incident
plane wave. For the case of lossless medium (i.e., $\I(k)=0$) above the conductor,
a local uniqueness result was obtained by Bao and Zhou in \cite{BZ} for the inverse problem
with one incident plane wave by establishing a lower bound of the first eigenvalue of
the $\curl\curl$ operator with the boundary condition (\ref{equation2}) in a bounded,
smooth convex domain in $\R^3.$ The stability of the inverse problem was also studied in \cite{BZ}.
Recently in \cite{BZZ09} it was proved that one incident plane wave is enough to uniquely
determine a bi-periodic polyhedral structure except for several extremely exceptional cases.
%
%If $\G_I=\emptyset$ and the surface $f(x_1,x_2)\in C^2(\R^2)$, Ammari \cite{AH}
%has proved that there is a unique quasi-periodic solution $E\in C^2(\Om)\cap(\overline{\Om})$
%and the solution depends continuously on $H^{-1/2}({\rm div},\G)$ boundary data.
%This generalizes the direct scattering results obtained in \cite{NS}.
%The uniqueness results of the inverse grating diffraction problems
%in \cite{Bao1994} and \cite{Kirsch1994} has also been generalized to
%the bi-periodic case by Ammari,H in \cite{AH}.
%Recently in \cite{BZZ09} it was proved that one incident plane wave is enough to uniquely
%determine a bi-periodic polyhedral structure except for several extremely exceptional cases.
%If $\G_I\neq\emptyset$ and $f(x_1,x_2)$ is  Lipschitz, the well-posedness of
%the direct problem was obtained in \cite{HQZ} with the aid of Hodge decomposition
%of $H(\curl,\Om_b)$. For the inverse problems, it is shown in \cite{HQZ} that,
%if the tangential components of scattered electric field on $\G_b$
%are identical for all incident electric dipoles $E^{in}(x;y)=\curl_x\curl_x\{PG(x,y)\}$
%with all $y\in\G_b$ and three linear independent vectors $P\in\R^3$,
%then their corresponding scattered surfaces $\G_j(j=1,2)$ and impedance
%coefficient $\lambda_j(j=1,2)$ on $\G_I$ must coincide. Here $G(x,y)$ is
%the free space quasi-periodic Green function (see Section \ref{sec2}).

Note that the inverse problem we are concerned in this paper involves in
the near field measurements since only a finite number of terms in (\ref{RE}) are
upwards propagating plane waves and the rest are evanescent modes that decay exponentially
with distance away from the grating. Thus we use near field data rather than far field data
to reconstruct the grating structure, which implies that the far field equation defined on the
unit ball of $\R^3$ for the non-periodic case must be replaced by a near field equation defined
on a plane above the structure. On the other hand, instead of using electromagnetic
Herglotz pairs in the case of bounded obstacle scattering problems,
we consider another kind of incident electric fields (see Section \ref{lemma} and
Remark \ref{rk1}) which lead to a denseness range result on the grating structure
since scattering occurs in a half space and the solution is $\al$-quasi-periodic
depending on the incident angle of the incident direction. This differs from the original
version of the linear sampling method which makes use of incident plane waves of all incident
directions with three linearly independent polarization directions (cf. \cite{CCM}).

The remaining part of the paper is organized as follows. Section \ref{sec2} is devoted to
the basic quasi-periodic function spaces used in the study of electromagnetic scattering
problems by periodic structures. Section \ref{lemma} gives several important lemmas
which are necessary for establishing the main result. The main result on the periodic
Linear Sampling Method and the numerical strategies on the implementation of the linear
sampling method are presented in Section \ref{lsm}.

\section{Basic Function Spaces}\label{sec2}
\setcounter{equation}{0}

In this section we introduce some quasi-periodic Sobolev spaces which are well-suited
for our problems. Due to the periodicity of the problem, the original problem can be
reduced to a problem in a single periodic cell of the grating profile.
To this end and for the subsequent analysis, we reformulate the following notations:
\ben
\G&=&\{x_3=f(x_1,x_2)\,|\,0<x_1,x_2<2\pi\},\\
\G_b&=&\{x_3=b\,|\,0<x_1,x_2<2\pi\},\\
\Om&=&\{x\in \R^3\,|\,x_3>f(x_1,x_2),0<x_1,x_2<2\pi\},\\
\Om_b&=&\{x\in\Om\,|\, x_3<b\},\\
\R^3_{\pi}&=&\{x\in\R^3: 0<x_1,x_2<2\pi\}
\enn
for any $b>\max\{f(x_1,x_2)\}$. We now introduce the scalar quasi-periodic Sobolev space:
\ben
H^1(\Om_b)=\{u(x)=\sum_{n\in\Z^2}u_n\exp[i(\al_n\cdot x+\beta_nx_3)]\,|\,
u\in L^2(\Om_b),\nabla u\in (L^2(\Om_b))^3,u_n\in\C\}.
\enn
Denote by $H^{\frac12}(\G_b)$ the trace space of $H^1(\Om_b)$ on $\G_b$ with the norm
\ben
||f||^2_{H^{\f{1}{2}}(\G_b)}=\sum_{n\in\Z^2}|f_n|^2(1+|\al_n|^2)^{\f{1}{2}},
  \qquad f\in H^{\frac12}(\G_b),
\enn
where $f_n=(f,\exp(i\al_n\cdot x))_{L^2(\G_b)}$.
Write $H^{-\frac12}(\G_b)=(H^{\frac{1}{2}}(\G_b))^\prime$, the dual space to $H^{\frac12}(\G_b)$.
We also need some vector spaces. Let
\ben
H(\curl,\Om_b)&=&\{E(x)=\sum_{n\in\Z^2}E_n\exp[i(\al_n\cdot x+\beta_nx_3)]\,|\,E_n\in\C^3,\\
              &&\quad\,E\in(L^2(\Om_b))^3,\,\curl E\in(L^2(\Om_b))^3\}
\enn
with the norm
\ben
||E||^2_{H(\curl,\Om_b)}=||E||^2_{L^2(\Om_b)}+||\curl E||^2_{L^2(\Om_b)}
\enn
and let
\ben
H_0(\curl,\Om_b)=\{E\in H(\curl,\Om_b),\,\nu\times E=0\,{\rm on }\,\G_b\}.
\enn
Define
\ben
X:=X(\Om_b,\G_I)=\{E\in H(\curl,\Om_b),\,\nu\times E|_{\G_I}\in L_t^2(\G_I)\}
\enn
with the norm
\ben
||E||^2_X=||E||^2_{H(\curl,\Om_b)}+||\nu\times E||^2_{L_t^2(\G_I)},
\enn
where $L_t^2(\G)=\{E\in(L^2(\G))^3,\,\nu\cdot E=0\,{\rm on}\,\G\}.$
For $x^\prime=(x_1,x_2,b)\in\G_b$, $s\in\R$ define
\ben
H_t^s(\G_b)&=&\{E(x^\prime)=\sum_{n\in\Z^2}E_n\exp(i\al_n\cdot x^\prime)\,|\,
              E_n\in\C^3,\,e_3\cdot E=0,\\
           &&\qquad \|E\|^2_{H^s(\G_b)}=\sum_{n\in\Z^2}(1+|\al_n|^2)^s|E_n|^2<+\infty\}\\
H_t^s(\dive,\G_b)&=&\{E(x^\prime)=\sum_{n\in\Z^2}E_n\exp(i\al_n\cdot x^\prime)\,|\,
                     E_n\in\C^3,\,e_3\cdot E=0,\\
                 &&\quad ||E||^2_{H^s(\dive,\G_b)}=\sum_{n\in\Z^2}(1+|\al_n|^2)^s
                        (|E_n|^2+|E_n\cdot\al_n|^2)<+\infty\}\\
H_t^s(\curl,\G_b)&=&\{E(x^\prime)=\sum_{n\in\Z^2}E_n\exp(i\al_n\cdot x^\prime)\,|\,
                     E_n\in\C^3,\,e_3\cdot E=0,\\
                 &&\quad ||E||^2_{H^s(\curl,\G_b)}=\sum_{n\in\Z^2}(1+|\al_n|^2)^s
                     (|E_n|^2+|E_n\times\al_n|^2)<+\infty\}
\enn
and write $L_t^2(\G_b)=H_t^0(\G_b).$ Recall that
\ben
H_t^{-{1}/{2}}(\dive,\G_b)=\{e_3\times E|_{\G_b}\,|\,E\in H(\curl,\Om_b)\}
\enn
and that the trace mapping from $H(\curl,\Om_b)$ to $H_t^{-{1}/{2}}(\dive,\G_b)$ is continuous and
surjective (see \cite{BCS} and the references there).
The trace space on the complementary part $\G_D $ of $X(\Om_b,\G_I)$ is
\ben
Y(\G_D)=\{f\in (H^{-{1}/{2}}(\G_D))^3\,|\,\exists E\in H_0(\curl,\Om_b)\,{\rm such}\,{\rm that}\,
           \nu\times E|_{\G_I}\in L_t^2(\G_I),\,\nu\times E|_{\G_D}=f\}
\enn
which is a Banach space with the norm
\ben
||f||^2_{ Y(\G_D)}&=&\inf\{||E||^2_{H(\curl,\Om_b)}+||\nu\times E||^2_{L_t^2(\G_I)}\,|\,\\
    &&\qquad E\in H_0(\curl,\Om_b),\,\nu\times E|_{\G_I}\in L_t^2(\G_I),\,\nu\times E|_{\G_D}=f\}.
\enn
An equivalent norm to $||\cdot||_{Y_{\G_D}}$ is given by (see \cite{CCM,CDZ,M})
\ben
|||f|||_1=\sup_{V\in X(\Om_b,\G_I)}\frac{|<f,V>_1|}{||V||_{X(\Om_b,\G_I)}},
\enn
where, for $E\in H_0(\curl,\Om_b)$ satisfying that $\nu\times E|_{\G_I}\in L_t^2(\G_I)$ and
$\nu\times E|_{\G_D}=f$, we have
\be\label{2.1}
<f,V>_1=\int_{\Om_b}(\curl E\cdot V-E\cdot\curl V)dx-\int_{\G_I}\nu\times E\cdot Vds,
     \quad\forall V\in X(\Om_b,\G_I).
\en
In particular, $Y(\G_D)$ is a Hilbert space, and (\ref{2.1}) can be considered as
a duality between $Y(\G_D)$ and its dual space $Y(\G_D)^\prime.$ From (\ref{2.1})
it can be seen that $\varphi\in Y(\G_D)^\prime$ can be extended as a function
$\widetilde{\varphi}\in H^{-1/2}_{\curl}(\G)$ defined on the whole boundary $\G$
such that $\widetilde{\varphi}|_{\G_I}\in L_t^2(\G_I)$.

We conclude this section with introducing the following free space Green function
$\Phi(x,y)$ for the Helmholtz equation $(\Delta+k^2)u=0$ in $\R^3$:
\ben
\Phi(x,y)=\frac{e^{ik|x-y|}}{4\pi|x-y|}
\enn
and the following free space $\al$-quasi-periodic Green function $G(x,y)$ for
the Helmholtz equation:
\be\label{Green}
G(x,y)=\frac{1}{8\pi^2}\sum_{n\in\mathbf{N}^2}\frac{1}{i\bt_n}
       \exp(i\al_n\cdot(x-y)+i\bt_n|x_3-y_3|)
\en
with $\al_n,\beta_n$ defined as in the introduction.

\section{Several Lemmas}\label{lemma}
\setcounter{equation}{0}

In this section we prove several important lemmas which are necessary for the proof of
the main theorem. We first define the incident electric field
\be\label{incident}
E^{in}(x;g):=\curl_x\curl_x\int_{\G_b}\overline{g(y)G(y,x)}ds(y)
\en
for $g\in L_t^2(\G_b)$. From the definition of $G(x,y)$ it is seen that $E^{in}(x;g)$
satisfies the radiation condition (\ref{RE}) in the region below $\G_b$.
This means that physically the above incident field propagates upward and does not
appear to be meaningful as incident waves.
Thus the total electric field corresponding to $E^{in}(x;g)$ can not be generated directly.
We will discuss how to solve the direct scattering problem for such incident waves
in the final section. For any $g\in L_t^2(\G_b)$ we next define a function
$(Hg)\in B:=Y(\G_D)\times L_t^2(\G_I)$ by
\ben
(Hg)(x)=\left\{\begin{array}{lll}
\nu(x)\times E^{in}(x;g),\quad\rm{on}\quad\G_D,\\
\nu(x)\times\curl E^{in}(x;g)-i\lambda E^{in}(x;g)_T,\quad\rm{on}\quad\G_I,\\
\end{array}\right.
\enn
where, for any vector field $V$, $V_T:=(\nu\times V)\times\nu$ denotes its tangential
component on a surface.

\begin{lemma}\label{denseness}
The range of $H$ is dense in $B$.
\end{lemma}

\begin{proof}
%We notice that both $E^{in}_{y}(x)(y\in\G_b)$ and $E^{in}_{z_0}(x)$
%propagate downward satifying Rayleigh Expansion (\ref{RE}) with
%$-\bt_n$ in $\R^3\backslash\overline{\Om}$. In view of the proof of
%Theorem \ref{th3.1}, the uniqueness and existence results with a
%negative impedance coefficient $\lambda$, still hold in the low
%region by introducing D-to-N map on an artificial boundary
%$\G_{-b}$ below the surface. For $y\in\G_b$ we next  define function
%$H_y^P(x)\in Y(\G_D)\times L_t^2(\G_I)$ by
%$$H_y^P(x)=\left \{\begin{array}{lll}
%\nu(x)\times E^{in}(x,y)\ \ \ \ \ \ \ \ \ {\rm on}\ \  \G_D,\\
%\nu(x)\times \curl_x E^{in}(x,y) +i\lambda E^{in}(x,y)_T\ \ \ \ {\rm on }\ \ \G_I.\\
%\end{array}\right. $$
%It is enough to prove that $\overline{{\rm Span}\{H^{P_i}_y: y\in
%\G_b, i=1,2,3\}}$ is dense in $Y(\G_D)\times L_t^2(\G_I)$.

For $f\times h\in B^{*}:=Y(\G_D)^{'}\times L_t^2(\G_I)$, we are going to prove that $f=0,h=0$
under the assumption that $<Hg,f\times h>_{B,B^{*}}=0$ for any $g\in L^{2}_{t}(\G_b)$.
Recalling that the duality between $Y(\G_D)$ and $Y(\G_D)^{'}$ is defined by (\ref{2.1})
and the duality between $L_t^2(\G_I)$ and $L_t^2(\G_I)$ is the $L^2$ scalar product,
we have
\ben
0&=&\int_{\G_D}\nu(x)\times\left[\curl_x\curl_x\int_{\G_b} g(y)G(y,x)ds(y)\right]\cdot f(x)ds(x)\\
 &&+\int_{\G_I}\left\{\nu(x)\times\left[\curl_x\curl_x\curl_x\int_{\G_b}g(y)G(y,x)ds(y)\right]\right.\\
 &&+i\lambda\left[\curl_x\curl_x\int_{\G_b}g(y)G(y,x)ds(y)\right]_T\left.\right\}\cdot h(x)ds(x)
\enn
Since $f\in Y(\G_D)^{'}$, there is an extension $\widetilde{f}\in H^{-1/2}_{\curl}(\G)$ of $f$,
defined on $\G$, satisfying that $\widetilde{f}|_{\G_I}\in L_t^2(\G_I)$.
Thus the above equation can be rewritten as
\ben
0&=&\int_{\G}\nu(x)\times\left[\curl_x\curl_x\int_{\G_b}g(y)G(y,x)ds(y)\right]
    \cdot\widetilde{f}(x)ds(x)\\
 &&-\int_{\G_I}\nu(x)\times\left[\curl_x\curl_x\int_{\G_b}g(y)G(y,x)ds(y)\right]
    \cdot\widetilde{f}(x)ds(x)\\
 &&+\int_{\G_I}\left\{\nu(x)\times\left[\curl_x\curl_x\curl_x\int_{\G_b}g(y)G(y,x)ds(y)\right]\right.\\
 &&+i\lambda\left[\curl_x\curl_x\int_{\G_b}g(y)G(y,x)ds(y)\right]_T\left.\right\}\cdot h(x)ds(x).
\enn
Making use of the vector identity:
\ben
\{\curl_x\curl_x[g(y)G(y,x)]\}\cdot h(x)=\{\curl_y\curl_y[h(x)G(y,x)]\}\cdot g(y),
\enn
we obtain by a direct computation that for any $g\in L^2_t(\G_b),$
\ben
     k^2\ds\int_{\G_b}E(y)\cdot g(y)ds(y)=0,
\enn
where, for $y\in\R^3_{\pi}\backslash\G$,
\ben
E(y)&=&\f{1}{k^2}\{\curl_y\curl_y\int_{\G}G(y,x)\widetilde{f}(x)\times\nu(x)ds(x)
        -\curl_y\curl_y\int_{\G_I}G(y,x)\widetilde{f}(x)\times\nu(x)ds(x)\\
   &&+k^2\curl_y\int_{\G_I}G(y,x)h(x)\times\nu(x)ds(x)
         +i\lambda\curl_y\curl_y\int_{\G_I}G(y,x)h(x)ds(x)\}.
\enn
Thus we have
\ben
\curl\curl E-k^2E=0&&\qquad y\in\R^3_{\pi}\backslash\G,\\
\nu\times E=0&&\qquad y\in\G_b.
\enn
It is clear that $E(y)$ propagates upward above $\G$ satisfying the Rayleigh expansion
radiation condition (\ref{RE}) and propagates downward below $\G$ satisfying the Rayleigh
expansion radiation condition (\ref{RE}) with $\al$ replaced by $-\al$.
By the uniqueness of the radiating solution to the exterior problem of the Maxwell equations
with the perfectly conducting condition and the analytic continuation of the solution of
the Maxwell equations, it follows that $E(y)\equiv 0$ for $y_3>f(y_1,y_2)$.
When $y\rightarrow\G$, the following jump relations hold on $\G$:
\ben
\nu\times E^{+}-\nu\times E^{-}&=&0,\qquad{\rm on}\quad\G_D,\\
i\lambda E^{+}-i\lambda E^{-}&=&-i\lambda h,\qquad{\rm on}\quad\G_I,\\
\nu\times\curl E^{+}-\nu\times\curl E^{-}&=&i\lambda h,\qquad{\rm on}\quad\G_I,
\enn
where the superscripts $+$ and $-$ indicate the limit obtained from $\Om$ and
$\R^3_{\pi}\backslash\overline{\Om}$, respectively.
It should be remarked that, since $\widetilde{f}\in H^{-1/2}(\curl,\G),$
the first integral over $\G$ in the definition of $E(y)$ is well defined
with $H^{-1/2}$-density (see \cite{Mw}) and the corresponding jump
conditions are interpreted in the sense of the $L^2$-limit.
Combining these jump relations and using the fact that
$\nu\times E^{+}=\nu\times\curl E^{+}=0$ lead to
\ben
\curl\curl E-k^2E=0&&\qquad y_3<f(y_1,y_2),\\
\nu\times E^{-}=0&&\qquad{\rm on}\quad\G_D,\\
\nu\times\curl E^{-}+i\lambda E^{-}_T=0&&\qquad{\rm on}\quad\G_I.
\enn
A similar argument as in \cite{HQZ} can be applied to the above problem to show
%show the existence of a unique solution to the above problem in the region below
%the doubly periodic structure with the impedance coefficient $-\lambda$ and
%the Dirichlet-to-Neumann map on the artificial boundary below the structure.
%Thus the application of uniqueness to the above equation yields
that $E(y)\equiv0$ for $y_3<f(x_1,x_2)$. Thus, we have
\ben
f =[\curl E]|_{\G_D}=0,\qquad\; h=-[\nu\times E]|_{\G_I}=0,
\enn
where $[\cdot]|_{\G_A}$ stands for the jump across $\G_A$ of a function with $A=D,I$.
The proof of Lemma \ref{denseness} is thus completed.
\end{proof}

The near field operator $N$ is defined by a bounded operator from $B$ into
$H^{-1/2}_t({\rm div},\G_b)$ which maps the boundary data $(h_1,h_2)\in B$ to
the tangential component $e_3\times E^s(x)|_{\G_b}$ of the near electric field.
Here, $E^s$ stands for the unique Rayleigh expansion radiating solution to the
Maxwell equations with the following boundary conditions:
\ben
\nu\times E^s=h_1\quad\mbox{on}\quad\G_D,\qquad\nu\times\curl E^s-i\lambda(E^s)_T=h_2
\quad\mbox{on}\quad\G_I.
\enn
By the well-posedness of the direct problem (see \cite{HQZ}) it is known
that $N$ is injective and bounded. Furthermore, $N$ is a compact operator.
To see this, we need the following periodic representation formula.

\begin{lemma}\label{representation}
Assume that $E$ satisfies the Rayleigh expansion radiation condition (\ref{RE})
and the Maxwell equations in $\Om$. Then for any $x\in\Om$ we have
\ben
E(x)=\curl_x\int_{\G}\nu(y)\times E(y)G(x,y)ds(y)
    +\frac{1}{k^2}\curl_x\curl_x\int_{\G}\nu(y)\times\curl E(y)G(x,y)ds(y),
\enn
where $G(x,y)$ is the quasi-periodic Green function defined by (\ref{Green}).
\end{lemma}

\begin{proof}
For arbitrarily fixed $x\in\Om$ and an arbitrary constant vector $P\in\R^3$
let $F(x,y)=PG(x,y)$ with $y\in\Om$.
Assume that $x\in\Om_b$ for some $b>0$. Denote by $B_{\delta}(x)$ the small ball
centered at $x$ with radius $\delta$ such that $B_{\delta}(x)\subset\Om_b$.
It is clear that both $E$ and $F(x,\cdot)$ satisfy the vector Helmhotz equation
in $\Om\backslash B_{\delta}(x)$. Using Green's second vector theorem and the
quasi-periodicity of $E$ and $F$ we have
\be\no
0&=&\int_{\Om_b\backslash B_{\delta}(x)}E(y)\cdot\triangle F(x,y)-\triangle E(y)\cdot F(x,y)dy\\ \no
&=&\left(-\int_{\G}+\int_{\G_b}+\int_{|y-x|=\delta}\right)\left\{\nu\times E\cdot\curl_yF+
    \nu\cdot E\rm{div}_yF-\nu\times F\cdot\curl E\right\}ds(y)\\ \label{a}
 &:=&-I_1+I_2+I_3.
\en
By a direct computation, we have
\be\no
I_1&:=&\int_{\G}\left\{\nu\times E\cdot\curl_y[PG(x,y)]
     +\nu\cdot E\rm{div}_y[PG(x,y)]-\nu\times[PG(x,y)]\cdot\curl E\right\}ds(y)\\ \no
   &=&\int_{\G}\left\{(\nu\times E)\times\nabla_yG(x,y)+(\nu\cdot E)\nabla_yG(x,y)
      +(\nu\times\curl E)G(x,y)\right\}ds(y)\cdot P\\ \label{b}
   &:=&\int_{\G}T_G(E)(y)ds(y)\cdot P
\en
with
\ben\no
T_G(E)(y)&=&(\nu\times E(y))\times\nabla_yG(x,y)+(\nu(y)\cdot E(y))\nabla_yG(x,y)
           +(\nu\times\curl E(y)) G(x,y)\\ \no
&=&\curl_x[\nu(y)\times E(y)G(x,y)]-\nabla_x[\nu(y)\cdot E(y)G(x,y)]
   +(\nu\times\curl E(y))G(x,y),
\enn
where we have used the fact that $\nabla_yG(x,y)=-\nabla_xG(x,y)$ to get the second equality,
\be\no
I_2&:=&\int_{\G_b}T_G(E)(y)ds(y)\cdot P\\\no
  &=&\int_{\G_b}\left\{[(\nabla_yG(x,y)\times P)\times e_3]\cdot E(y)
      +(\nabla_yG(x,y)\cdot P)(e_3\cdot E(y))\right.\\ \label{b+}
  &&-PG(x,y)\cdot(\curl E(y)\times e_3)\left.\right\}ds(y)=0,
\en
where the last equality follows from the Rayleigh expansion condition (\ref{RE}) for $E(y)$,
the definition of $G(x,y)$ and the fact that $\rm{div} E\equiv0$ in $\Om$,
\be\no
I_3&:=&\int_{|y-x|=\delta}T_G(E)(y)ds(y)\cdot P\\ \label{c}
   &=&\left(\int_{|y-x|=\delta}T_{G-\Phi}(E)(y)ds(y)
     +\int_{|y-x|=\delta}T_{\Phi}(E)(y)ds(y)\right)\cdot P
\en
Since $G(x,y)-\Phi(x,y)$ is a $C^{\infty}$-function with respect to $y$
in $B_{\delta}(x)$ (see \cite{NS}), we have
$\int_{|y-x|=\delta}T_{G-\Phi}(E)(y)ds(y)\rightarrow 0$ as $\delta\rightarrow0$.
The application of the mean value theorem yields that
$\int_{|y-x|=\delta}T_{\Phi}(E)(y)ds(y))\rightarrow E(x)$ as $\delta\rightarrow 0$ (cf. \cite{CK}).
Thus it follows from (\ref{a})-(\ref{c}) that
\ben
E(x)&=&\curl_x\int_{\G}\nu(y)\times E(y)G(x,y)ds(y)+\int_{\G}\nu(y)\times\curl E(y)G(x,y)ds(y)\\
     &&-\nabla_x\int_{\G}\nu(y)\cdot E(y)G(x,y)ds(y)
\enn
which is analogous to the well-know non-periodic Stratton-Chu representation
theorem (\cite[Theorem 6.1]{CK}). Finally, the application of the Stokes theorem
together with the vector identity $\curl\curl=-\triangle+\nabla(\nabla\cdot)$ gives
the desired result.
\end{proof}

It is seen from Lemma \ref{representation} and the well-posedness of the
direct scattering problem that $N$ is a composition of a bounded operator
mapping the boundary data into the scattered field with a compact operator
taking the scattered field to its tangential component of $E^s$ on $\G_b$.
Thus $N$ is compact. We now prove, with the help of Lemma \ref{representation},
that $N$ has a dense range in $H^{-1/2}_t({\rm div},\G_b)$.

\begin{lemma}\label{denseness2}
The set $\{N(\varphi,\psi)\;|\;\varphi\in Y(\G_D),\psi\in L_t^2(\G_I)\}$ is dense
in $H^{-1/2}_t({\rm div},\G_b).$
\end{lemma}

\begin{proof}
Let $h\in H^{-1/2}_t(\curl,\G_b)=H^{-1/2}_t({\rm div},\G_b)^\prime$ satisfy
\be\label{adjoint}
<N(\varphi,\psi),h>=0\qquad\forall\,\varphi\in Y(\G_D),\,\psi\in L_t^2(\G_I),
\en
where $<\cdot,\cdot>$ denotes the duality between $H^{-1/2}_t(\curl,\G_b)$
and $H^{-1/2}_t({\rm div},\G_b)$. Then it is sufficient to prove that $h=0$.
By the definition of N and the well-posedness of the direct scattering problem
there exists a unique $E\in H_{loc}(\curl,\Om)$ satisfying the Rayleigh expansion
radiation condition (\ref{RE}) such that $N(\varphi,\psi)=e_3\times E$ on $\G_b$.
From Lemma \ref{representation} it follows that
\ben
<N(\varphi,\psi),h>&=&\int_{\G_b}e_3\times E(x)\cdot\overline{h(x)}ds(x)\\
 &=&\int_{\G_b}(\overline{h(x)}\times e_3)\cdot\left\{\curl_x
    \int_{\G}\nu(y)\times E(y)G(x,y)ds(y)\right\}ds(x)\\
 &&+\frac{1}{k^2}\int_{\G_b}(\overline{h(x)}\times e_3)\cdot
   \left\{\curl_x\curl_x\int_{\G}\nu(y)\times\curl E(y)G(x,y)ds(y)\right\}ds(x)\\
 &:=&I_1+I_2.
\enn
Interchanging the order of integration gives
\ben
I_1&=&-\int_{\G}\nu(y)\times E(y)\cdot\left\{\curl_y
      \int_{\G_b}G(x,y)\overline{h(x)}\times e_3ds(x)\right\}ds(y),\\
I_2&=&\frac{1}{k^2}\int_{\G}\nu(y)\times\curl E(y)\cdot\left\{\curl_y\curl_y
      \int_{\G_b}G(x,y)\overline{h(x)}\times e_3ds(x)\right\}ds(y).
\enn
Let
\ben
F(y):=\frac{1}{k^2}\curl_y\curl_y\int_{\G_b}G(x,y)\overline{h(x)}\times e_3ds(x),
\qquad y\in\R^3_{\pi}\backslash\G_b.
\enn
Then, since $\curl\curl=-\triangle+\nabla(\nabla\cdot)$, we have
\ben
\curl F(y)= -\curl_y\int_{\G_b}\overline{h(x)}\times e_3G(x,y)ds(x).
\enn
Thus
\be\label{e}
<N(\varphi\times\psi),h>=\int_{\G}\nu(y)\times E(y)\cdot\curl F(y)
-\nu(y)\times F(y)\cdot\curl E(y)ds(y).
\en
Let $\tilde{E}$ be the $-\al$-quasi-periodic Rayleigh expansion
radiating solution to the problem:
\be\no
\curl\curl\tilde{E}-k^2\tilde{E}&=&0
\qquad\quad\quad\quad\quad\quad\quad\quad\rm{in}\quad\Om,\\ \label{E1}
\nu\times\tilde{E}&=&\nu\times F\qquad\quad\quad\quad\quad\quad\rm{on}\quad\G_D,\\ \label{E2}
\nu\times\curl\tilde{E}-i\lambda\tilde{E}_T&=&\nu\times\curl F-i\lambda F_T\quad\rm{on}\quad\G_I.
\en
From Green's second vector theorem and the Rayleigh expansion of $\tilde{E}$ and $E$
it follows that
\ben
&&\int_{\G}\nu(y)\times E(y)\cdot\curl\tilde{E}(y)
  -\nu(y)\times\tilde{E}(y)\cdot\curl E(y)ds(y)\\ \no
&&\qquad=\int_{\G_b}\nu(y)\times E(y)\cdot\curl\tilde{E}(y)
  -\nu(y)\times \tilde{E}(y)\cdot\curl E(y)ds(y)=0,
\enn
which, in conjunction with the boundary conditions (\ref{E1}) and (\ref{E2}),
$\nu\times E=\varphi$ on $\G_D$ and $\nu\times\curl E-i\lambda E_T=\psi$ on $\G_I$,
implies that
\ben
\int_{\G_D}\nu\times F\cdot\curl E+\int_{\G_I}(\nu\times\curl F-i\lambda F_T)\cdot E
=\int_{\G_D}\varphi\cdot\curl\tilde{E}+\int_{\G_I}\psi\cdot\tilde{E}.
\enn
This together with (\ref{e}) yields
\ben
<N(\varphi,\psi),h>&=&\int_{\G_D}\varphi\cdot\curl F-\nu\times F\cdot\curl E
 +\int_{\G_I}\psi\cdot F- (\nu\times\curl F-i\lambda F_T)\cdot E\\
&=&\int_{\G_D}\varphi\cdot\left[\curl F-\curl\tilde{E}\right]
 +\int_{\G_I}\psi\cdot(F-\tilde{E}).
\enn
It is seen from the above identity that the conjugate operator of $N$ is given by
\ben
N^*h=\left((\overline{\curl F-\curl \tilde{E}})_T, (\overline{F-\tilde{E}})_T\right)\in B.
\enn
Combining (\ref{adjoint}) and the boundary conditions (\ref{E1}) and (\ref{E2}) gives
$\nu\times F=\nu\times\tilde{E}$ and $\nu\times\curl F=\nu\times\curl\tilde{E}$ on $\G$.
By Holmgren's uniqueness theorem, $F\equiv\tilde{E}$ in $\Om_b$.
Since $\nu\times F=\nu\times\tilde{E}$ on $\G_b$ and both $F$ and $\tilde{E}$ satisfy
the $-\al$-quasi-periodic Rayleigh expansion radiation condition for $x_3>b$,
it follows from the uniqueness result for the exterior Dirichlet problem that
$F\equiv\tilde{E}$ for $x_3>b$. Now, in view of the fact that $\tilde{E}$ is analytic
in $\Om$, we have by the jump relation of $\curl F(y)$ as $y\rightarrow\G_b$ that
\ben
\overline{h}=\left[\curl F^+-\curl F^-\right]|_{\G_b}
=\left[\curl E^+-\curl E^-\right]|_{\G_b}=0,
\enn
which completes the proof of the lemma.
\end{proof}

\section{The Linear Sampling Method}\label{lsm}
\setcounter{equation}{0}

For $g\in L_t^2(\G_b)$ consider $E^{in}(x;g)$ defined by (\ref{incident}) as
incident waves. Denote by $E^s(x;g)$ the scattered solution of the problem
(\ref{equation1})-(\ref{equation4}) corresponding to $E^{in}(x;g)$.
To derive a periodic version of the linear sampling method consider
the following near field equation:
\be\label{nearfield}
\mathcal{F}(g_z):=\int_{\G_b}e_3\times E^s(x,g_z)ds(x)
=e_3\times\curl_x\curl_x\{PG(x,z)\}\qquad{\rm on}\quad\G_b,
\en
where $z\in\{z\in\R^3|0<z_3<b\}$ and $P\in\R^3$ is a polarization vector.
It is clear that
\be\label{decomposition}
NH(g)=-\mathcal{F}(g).
\en

\begin{theorem}\label{PLSM}
Assume that $\G$ is Lipschitz continuous with the dissection $\G=\G_D\cup\Sigma\cup\G_I$
and $\G_I\neq\emptyset$.

(1) If $z\in\R^3_{\pi}\backslash\overline{\Om}$, then for any $\epsilon>0$
there exists a $g_{z,P}^{\epsilon}\in L_{t}^{2}(\G_b)$ such that
\ben
||(\mathcal{F}g_{z,P}^{\epsilon})-\nu\times\curl\curl\{PG(\cdot,z)\}||_{L^2_t(\G_b)}<\epsilon
\enn
and
\ben
||g_{z,P}^{\epsilon}||_{L_{t}^{2}(\G_b)}\rightarrow\infty\qquad{\rm as}\quad z\rightarrow\G^-.
\enn

(2) If $z\in\Om$, then for any $\epsilon>0$ and $\delta>0$ there exists a
$g_{z,P}^{\epsilon,\delta}\in L_{t}^{2}(\G_b)$ such that
\ben
||(\mathcal{F}g_{z,P}^{\epsilon,\delta})
-\nu\times\curl\curl\{PG(\cdot,z)\}||_{L^{2}_{t}(\G_b)}<\epsilon+\delta
\enn
and
\ben
||g_{z,P}^{\epsilon,\delta}||_{L_{t}^{2}(\G_b)}\rightarrow\infty\qquad{\rm as}\quad\delta\rightarrow0.
\enn
\end{theorem}

\begin{proof}
(1) Let $z\in \R^3_{\pi}\backslash \overline{\Om}$. In this case,
$e_3\times\curl_x\curl_x\{PG(x,z)\}|_{\G_b}$ is in the range of $N$
since it is the tangential component of the electric field $E^{s}_{P,z}:=\curl_x\curl_x\{PG(x,z)\}$
which is a solution of the exterior mixed boundary value problem with boundary data
$h_1=\nu\times E^{s}_{P,z}$ on $\G_D$ and $h_2=\nu\times\curl E^{s}_{P,z}-i\lambda(E^{s}_{P,z})_T$
on $\G_I$, that is,
\be\label{boundary}
e_3\times\curl\curl\{PG(x,z)\}|_{\G_b}=N(h_1,h_2).
\en
It can then be seen from the denseness of the range of $H$ that, for every
$\epsilon$>0 there is a $g^{\epsilon}_{P,z}:=g(\cdot;\epsilon,P,z)\in L_t^2(\G_b)$ such that
\be\label{H}
||H(g^{\epsilon}_{P,z})+(h_1,h_2)||_{Y(\G_D)\times L_t^2(\G_I)}<\epsilon.
\en
The boundedness of $N$ implies that
\ben
||NHg^{\epsilon}_{P,z}+N(h_1,h_2)||_{L_t^2(\G_b)}<C\epsilon
\enn
for some positive constant $C$. From this, (\ref{decomposition}) and (\ref{boundary})
it follows that
\ben
||(\mathcal{F}g_{z,P}^{\epsilon})-\nu\times\curl\curl\{PG(\cdot,z)\}||_{L^{2}_{t}(\G_b)}<C\epsilon.
\enn
Furthermore, if $z\rightarrow\G^{-}$, then we have
\ben
|| (h_1,h_2)||_{Y(\G_D)\times L_t^2(\G_I)}\rightarrow\infty
\enn
due to the singularity of $h_1$ and $h_2$ as $z\rightarrow\G^{-}$.
This, together with (\ref{H}), gives rise to
\ben
\lim_{z\rightarrow\G^{-}}||Hg^{\epsilon}_{P,z} ||_{Y(\G_D)\times L_t^2(\G_I)}=\infty,
\enn
which together with the boundedness of $H$ implies that
\ben
\lim_{z\rightarrow\G^{-}}||g^{\epsilon}_{P,z}||_{L_t^2(\G_b)}=\infty.
\enn

(2) Let $z\in\Om$. In this case, $e_3\times\curl\curl\{PG(x,z)\}|_{\G_b}$ is not in
the range of $N$ since, otherwise, $\curl\curl\{PG(x,z)\}$ will be a solution to
the Maxwell equations in $\R^3_{\pi}\backslash\overline{\Om}$ which is impossible
due to its singularity at $z$. However, using the Tikhonov regularization, we can
construct a regularized solution to the near field equation (\ref{nearfield}) since,
by Lemmas \ref{denseness} and \ref{denseness2}, $\mathcal{F}$ is compact and has a dense range.
Specifically, for an arbitrary $\delta>0$ there exist functions
$(h_{1,P,z}^{\delta},h_{2,P,z}^{\delta})\in Y(\G_D)\times L_t^2(\G_I)$
corresponding to some parameter $\alpha=\alpha(\delta)$ chosen by a regularization
strategy (e.g., the Morozov discrepancy principle) such that
\be\label{inequa1}
||N(h_{1,P,z}^{\delta},h_{2,P,z}^{\delta})
-\nu\times\curl\curl\{PG(\cdot,z)\}||_{L^{2}_{t}(\G_b)}<\delta.
\en
Furthermore, using the regularization strategy and the Picard theorem (see \cite{CK})
we get
\be\label{regularization}
\lim_{\delta\rightarrow 0}\alpha(\delta)=0,\qquad\quad
\lim_{\delta\rightarrow 0}||(h_{1,P,z}^{\delta},h_{2,P,z}^{\delta})||_{Y(\G_D)
\times L_t^2(\G_I)}=\infty.
\en
Then by Lemma \ref{denseness} and the boundedness of $N$, for any $\epsilon>0$
it is possible to find a $g_{P,z}^{\epsilon,\delta}\in L_t^2(\G_b)$ such that
\be\label{inequa2}
||N(Hg_{P,z}^{\epsilon,\delta})-N(h_{1,P,z}^{\delta},h_{2,P,z}^{\delta})||_{L_t^2(\G_b)}<\epsilon.
\en
Thus we have from (\ref{inequa1}) and (\ref{inequa2}) that
\ben
||(\mathcal{F}g_{P,z}^{\epsilon,\delta})
-\nu\times\curl\curl\{PG(\cdot,z)\}||_{L^{2}_{t}(\G_b)}<\epsilon+\delta.
\enn
Finally, by (\ref{regularization}) and (\ref{inequa2}) in conjunction with
the boundedness of $H$ and $N$, we have
\ben
||g_{z,P}^{\epsilon,\delta}||_{L_{t}^{2}(\G_b)}\rightarrow\infty
\qquad{\rm as}\quad\delta\rightarrow0.
\enn
The proof is thus completed.
\end{proof}

We now discuss some numerical strategies on the implementation of the
above linear sampling method.

As stated in Section \ref{lemma}, the incident waves $E^{in}(x;g)$
defined by (\ref{incident}) are not of physical relevance since they
propagate away from the surface. Thus $E^s(x;g(y))$,
the scattered field corresponding to $E^{in}_g(x)$, can not be generated directly.
In what follows, we make use of the method of Arens and Kirsch \cite{AK}
to generate $E^s(x;g)$.  We first examine that
\be\no
G(x,y)-\overline{G(y,x)}&=&\frac{1}{8\pi^2}
  \left\{\sum_{\al_n\leq k}\frac{1}{i\bt_n}e^{i(\al_n\cdot(x-y)-\bt_n(y_3-x_3))}
  +\sum_{\al_n\leq k}\frac{1}{i\bt_n}e^{i(\al_n\cdot(x-y)+\bt_n(y_3-x_3))}\right\}\\
&:=&\Delta^{(U)}(x;y)+\Delta^{(D)}(x;y)\label{difference}
\en
for $y\in\G_b$ and $x\in\Om_b$.
Note that $\Delta^{(U)}(x;y)$ and $\Delta^{(D)}(x;y)$ are upward and downward
propagating modes respectively. Set
\ben
E^{(U)}(x;g)&:=&\curl_x\curl_x\int_{\G_b}\overline{g(y)}\Delta^{(U)}(x;y)ds(y),\\
E^{(D)}(x;g)&:=&\curl_x\curl_x\int_{\G_b}\overline{g(y)}\Delta^{(D)}(x;y)ds(y),\\ \label{d}
\widetilde{E}^{in}(x;g)&:=&\curl_x\curl_x\int_{\G_b}\overline{g(y)}
 \left\{G(x;y)-\Delta^{(D)}(x,y)\right\}ds(y).
\enn
Clearly, $\widetilde{E}^{in}(x;g)$ is propagating towards the scattering surface,
so the corresponding unique scattered filed $\widetilde{E}^s(x;g)$ can be computed directly.
It is seen from (\ref{difference}) and the boundary value of $\widetilde{E}^s(x,g)$ that
\ben
\nu\times\left\{\widetilde{E}^s(x;g)|_{\G_D}+E^{(U)}(x;g)|_{\G_D}\right\}
&=&\nu\times\left\{\widetilde{E}^s(x;g)|_{\G_D}+\widetilde{E}^{in}(x;g)|_{\G_D}
   -E^{in}(x;g)|_{\G_D}\right\}\\
&=&-\nu\times E^{in}(x;g)|_{\G_D}.
\enn
Similarly, we have
\ben
&&\nu\times\curl\left[\widetilde{E}^s(x;g)+E^{(U)}(x;g)\right]
  -i\lambda\left[\widetilde{E}^s(x;g)+E^{(U)}(x;g)\right]_T\\
&&\quad=\nu\times\curl\left[\widetilde{E}^s(x;g)+\widetilde{E}^{in}(x;g)-E^{in}(x;g)\right]
  -i\lambda\left[\widetilde{E}^s(x;g)+\widetilde{E}^{in}(x;g)-E^{in}(x;g)\right]_T\\
&&\quad=-\left\{\nu\times\curl E^{in}(x;g)-i\lambda E^{in}(x;g)_T\right\}.
\enn
It follows from the uniqueness of the direct scattering problem that
$E^s(x;g)=\widetilde{E}^s(x;g)+E^{(U)}(x;g).$
Thus we can exactly generate $E^s(x,g)$ using the incident field $\widetilde{E}^{in}(x;g)$.

\begin{remark}\label{rk1} {\rm
The reason why we use $E^{in}(x;g)$ defined by (\ref{incident}) as incident waves
to reconstruct the periodic structure is that, by Lemma \ref{denseness}, such kinds of
incident waves lead to a dense range in $B$ of the operator $H$.
It should be remarked that, if $\G_I=\emptyset,$ that is, the total electric field
$E(x)$ satisfies the perfectly conducting boundary condition $\nu\times E=0$ on $\G$,
then we are allowed to choose the following field as incident waves:
\ben
E^{in}_d(x;g):=\curl_x\curl_x\int_{\G_b}g(y)G(x,y)ds(y),\qquad x_3<b
\enn
with $g\in L_t^2(\G_b)$. This kind of incident fields leads to a dense range
in $H^{-1/2}({\rm div},\G)$ of the operator mapping $g$ into the tangential
component on $\G$ of $E^{in}_d(x;g)$.
Since $E^{in}_d(x;g)$ propagates downward in $\Om_b$, the corresponding scattered field
can be produced directly. Thus the above strategy of generating the scattered
field corresponding to $E^{in}(x;g)$ can be avoided.
}
\end{remark}

Our reconstruction algorithm consists of the following three steps:
\begin{description}
\item {\bf Step 1.} Select a mesh of sampling points in a computing region
$\Sigma_b=\{x\in\R^3_\pi|0<x_3<b,0<x_1,x_2<2\pi\}$ which contains the grating surface.
\item {\bf Step 2.} Making use of the Tikhonov regularization and the Morozov discrepancy
principle to compute an approximate solution $g_{z,P}^{\epsilon}$ to the near field
equation (\ref{nearfield}).
\item {\bf Step 3.} Consider $||g_{z,P}^{\epsilon}||_{L_{t}^{2}(\G_b)}$ as an indicator
function of the sampling points $z$ and get the contour plot of
$||g_{z,P}^{\epsilon}||_{L_{t}^{2}(\G_b)}$ as a function of $z$.
\end{description}

\begin{remark}\label{rk2} {\rm
The above algorithm has been implemented in \cite{HQZ2} for the two-dimensional
TE polarization case, where the Maxwell equations are replaced by the scalar Helmholtz
equation and the boundary conditions on $\G_D$ and $\G_I$ are replaced with the
Dirichlet and impedance conditions, respectively. The numerical reconstruction results
presented in \cite{HQZ2} have shown the efficiency of the algorithm.
The implementation of the above algorithm for the three-dimensional case of the
full Maxwell equations is still in progress.
}
\end{remark}

\begin{remark}\label{rk3} {\rm
It follows form Theorem \ref{PLSM} that the indicator function can be used to
characterize the different region below and above the grating surface.
The numerical implementation of the linear sampling method can be found in \cite{CCM}
for inverse electromagnetic scattering problems by general bounded obstacles,
which has been proven to be very successful and effective once the necessary direct
scattering data are available. In order to get a better reconstruction result
the mesh in Step $1$ must be fine so that the characterization of the grating
surface could be clear. But this also increases the computational cost since
the near field equation (\ref{nearfield}) must be solved at each sampling point.
To avoid this, a multilevel linear sampling method was proposed in \cite{LLZ}
for the inverse acoustic obstacle scattering problems.
%studying the choice of the mesh points in Step $1$
The multilevel linear sampling method has been shown to be effective and
to possess asymptotically optimal computational complexity and thus provides
a fast numerical technique to implement the linear sampling method.
}
\end{remark}

%\section{Apendix}
%\setcounter{equation}{0}

%The following Picard Theorem can be seen in \cite{K} or \cite{CK}.
%\begin{theorem}\label{Picard} (Picard)
%Let $K:X\rightarrow Y$ be a linear compact operator with singular
%system $(\mu_j,x_j,y_j)$, then the equation
%\be\label{5.2}
%Kx=y
%\en
%is solvable if and only if
%\ben
%y\in\mathcal{N}(K^*)^{\bot}\ {\rm and}\ \displaystyle\sum_{j}\f{1}{\mu_j^2}|(y,y_j)|^2<\infty
%\enn
%where $\mathcal{N}(A)$ denotes the null space of $A$. In this case
%\ben
%x=\displaystyle\sum_{j}\f{1}{\mu_j^2}(y,y_j)x_j
%\enn
%is a solution of (\ref{5.2}).
%\end{theorem}

\section*{Acknowledgements}

This work was supported by the NNSF of China under grant No. 10671201.

\end{document}